\font\fifteenrm=cmr10 scaled\magstep2 
\font\fifteeni=cmmi10 scaled\magstep2
\font\fifteensy=cmsy10 scaled\magstep2
\font\fifteenbf=cmbx10 scaled\magstep2
\font\fifteentt=cmtt10 scaled\magstep2
\font\fifteenit=cmti10 scaled\magstep2
\font\fifteensl=cmsl10 scaled\magstep2
\font\fifteenam=msam10 scaled\magstep2
\font\fifteenbm=msbm10 scaled\magstep2
\font\fifteenex=cmex10 scaled\magstep2
\font\fifteensc=cmcsc10 scaled\magstep2 
\font\twelverm=cmr10 at 12pt
\font\twelvei=cmmi10 at 12pt
\font\twelvesy=cmsy10 at 12pt
\font\twelvebf=cmbx10 at 12pt
\font\twelvett=cmtt10 at 12pt
\font\twelveit=cmti10 at 12pt
\font\twelvesl=cmsl10 at 12pt
\font\twelveam=msam10 at 12pt
\font\twelvebm=msbm10 at 12pt
\font\twelveex=cmex10 at 12pt
\font\twelvesc=cmcsc10 at 12pt
\font\elevenrm=cmr10 scaled\magstephalf 
\font\eleveni=cmmi10 scaled\magstephalf
\font\elevensy=cmsy10 scaled\magstephalf
\font\elevenbf=cmbx10 scaled\magstephalf
\font\eleventt=cmtt10 scaled\magstephalf
\font\elevenit=cmti10 scaled\magstephalf
\font\elevensl=cmsl10 scaled\magstephalf
\font\elevenam=msam10 scaled\magstephalf
\font\elevenbm=msbm10 scaled\magstephalf
\font\elevenex=cmex10 scaled\magstephalf
\font\elevensc=cmcsc10 scaled\magstephalf
\font\tenrm=cmr10
\font\teni=cmmi10
\font\tensy=cmsy10
\font\tenbf=cmbx10
\font\tentt=cmtt10
\font\tenit=cmti10
\font\tensl=cmsl10
\font\tenam=msam10
\font\tenbm=msbm10
\font\tenex=cmex10
\font\tensc=cmcsc10
\font\ninerm=cmr9
\font\ninei=cmmi9
\font\ninesy=cmsy9
\font\ninebf=cmbx9
\font\ninett=cmtt9
\font\nineit=cmti9
\font\ninesl=cmsl9
\font\nineam=msam9
\font\ninebm=msbm9
\font\nineex=cmex9
\font\ninesc=cmcsc9
\font\eightrm=cmr8
\font\eighti=cmmi8
\font\eightsy=cmsy8
\font\eightbf=cmbx8
\font\eighttt=cmtt8
\font\eightit=cmti8
\font\eightsl=cmsl8
\font\eightam=msam8
\font\eightbm=msbm8
\font\eightex=cmex8
\font\eightsc=cmcsc8
\font\sevenrm=cmr7
\font\seveni=cmmi7
\font\sevensy=cmsy7
\font\sevenbf=cmbx7

\font\sevenam=msam7
\font\sevenbm=msbm7

\font\sixrm=cmr6
\font\sixi=cmmi6
\font\sixsy=cmsy6

\font\sixam=msam6
\font\sixbm=msbm6

\font\fiverm=cmr5
\font\fivei=cmmi5
\font\fivesy=cmsy5

\font\fiveam=msam5
\font\fivebm=msbm5

\font\fourrm=cmr5 at 4pt
\font\fouri=cmmi5 at 4pt
\font\foursy=cmsy5 at 4pt

\font\fouram=msam5 at 4pt
\font\fourbm=msbm5 at 4pt

\skewchar\twelvei='177 \skewchar\eleveni='177\skewchar\teni='177
\skewchar\ninei='177 \skewchar\eighti='177\skewchar\seveni='177 
\skewchar\sixi='177 \skewchar\fivei='177 \skewchar\fouri='177
\skewchar\twelvesy='60 \skewchar\elevensy='60 \skewchar\tensy='60
\skewchar\ninesy='60 \skewchar\eightsy='60 \skewchar\sevensy='60 
\skewchar\sixsy='60 \skewchar\fivesy='60 \skewchar\foursy='60
\newfam\itfam
\newfam\slfam
\newfam\bffam
\newfam\ttfam
\newfam\scfam
\newfam\amfam
\newfam\bmfam
\def\eightbig#1{{\hbox{$\left#1\vbox to 6.5pt{}\voidright $}}}
\def\eightBig#1{{\hbox{$\left#1\vbox to 7.5pt{}\voidright $}}}
\def\eightbigg#1{{\hbox{$\left#1\vbox to 10pt{}\voidright $}}}
\def\eightBigg#1{{\hbox{$\left#1\vbox to 13pt{}\voidright $}}}
\def\ninebig#1{{\hbox{$\left#1\vbox to 7.5pt{}\voidright $}}}
\def\nineBig#1{{\hbox{$\left#1\vbox to 8.5pt{}\voidright $}}}
\def\ninebigg#1{{\hbox{$\left#1\vbox to 11.5pt{}\voidright $}}}
\def\nineBigg#1{{\hbox{$\left#1\vbox to 14.5pt{}\voidright $}}}
\def\tenbig#1{{\hbox{$\left#1\vbox to 8.5pt{}\voidright $}}}
\def\tenBig#1{{\hbox{$\left#1\vbox to 9.5pt{}\voidright $}}}
\def\tenbigg#1{{\hbox{$\left#1\vbox to 12.5pt{}\voidright $}}}
\def\tenBigg#1{{\hbox{$\left#1\vbox to 16pt{}\voidright $}}}
\def\elevenbig#1{{\hbox{$\left#1\vbox to 9pt{}\voidright $}}}
\def\elevenBig#1{{\hbox{$\left#1\vbox to 10.5pt{}\voidright $}}}
\def\elevenbigg#1{{\hbox{$\left#1\vbox to 14pt{}\voidright $}}}
\def\elevenBigg#1{{\hbox{$\left#1\vbox to 17.5pt{}\voidright $}}}
\def\twelvebig#1{{\hbox{$\left#1\vbox to 10pt{}\voidright $}}}
\def\twelveBig#1{{\hbox{$\left#1\vbox to 11pt{}\voidright $}}}
\def\twelvebigg#1{{\hbox{$\left#1\vbox to 15pt{}\voidright $}}}
\def\twelveBigg#1{{\hbox{$\left#1\vbox to 19pt{}\voidright $}}}
\def\fifteenbig#1{{\hbox{$\left#1\vbox to 12pt{}\voidright $}}}
\def\fifteenBig#1{{\hbox{$\left#1\vbox to 13.5pt{}\voidright $}}}
\def\fifteenbigg#1{{\hbox{$\left#1\vbox to 18pt{}\voidright $}}}
\def\fifteenBigg#1{{\hbox{$\left#1\vbox to 23pt{}\voidright $}}}
\def\voidright{\right.\nulldelimiterspace=0pt \mathsurround=0pt }
\def\fifteenpoint{
  \textfont0=\fifteenrm \scriptfont0=\twelverm \scriptscriptfont0=\tenrm
  \def\rm{\fam0 \fifteenrm}%
  \textfont1=\fifteeni \scriptfont1=\twelvei \scriptscriptfont1=\teni
  \textfont2=\fifteensy \scriptfont2=\twelvesy \scriptscriptfont2=\tensy
  \textfont3=\fifteenex \scriptfont3=\fifteenex \scriptscriptfont3=\fifteenex
  \def\it{\fam\itfam\fifteenit}\textfont\itfam=\fifteenit
  \def\sl{\fam\slfam\fifteensl}\textfont\slfam=\fifteensl
  \def\bf{\fam\bffam\fifteenbf}\textfont\bffam=\fifteenbf 
    \scriptfont\bffam=\twelvebf\scriptscriptfont\bffam=\tenbf
  \def\tt{\fam\ttfam\fifteentt}\textfont\ttfam=\fifteentt
  \def\sc{\fam\scfam\fifteensc}\textfont\scfam=\fifteensc
  \def\am{\fam\amfam\fifteenam}\textfont\amfam=\fifteenam
    \scriptfont\amfam=\twelveam\scriptscriptfont\amfam=\tenam
  \def\bm{\fam\bmfam\fifteenbm}\textfont\bmfam=\fifteenbm
    \scriptfont\bmfam=\twelvebm\scriptscriptfont\bmfam=\tenbm
  \baselineskip=21pt \rm
  \let\big=\fifteenbig\let\Big=\fifteenBig\let\bigg=\fifteenbigg
  \let\Bigg=\fifteenBigg}
\def\twelvepoint{
  \textfont0=\twelverm \scriptfont0=\ninerm \scriptscriptfont0=\sevenrm
  \def\rm{\fam0 \twelverm}%
  \textfont1=\twelvei \scriptfont1=\ninei \scriptscriptfont1=\seveni
  \textfont2=\twelvesy \scriptfont2=\ninesy \scriptscriptfont2=\sevensy
  \textfont3=\twelveex \scriptfont3=\twelveex \scriptscriptfont3=\twelveex
  \def\it{\fam\itfam\twelveit}\textfont\itfam=\twelveit
  \def\sl{\fam\slfam\twelvesl}\textfont\slfam=\twelvesl
  \def\bf{\fam\bffam\twelvebf}\textfont\bffam=\twelvebf 
    \scriptfont\bffam=\ninebf\scriptscriptfont\bffam=\sevenbf
  \def\tt{\fam\ttfam\twelvett}\textfont\ttfam=\twelvett
  \def\sc{\fam\scfam\twelvesc}\textfont\scfam=\twelvesc
  \def\am{\fam\amfam\twelveam}\textfont\amfam=\twelveam
    \scriptfont\amfam=\nineam\scriptscriptfont\amfam=\sevenam
  \def\bm{\fam\bmfam\twelvebm}\textfont\bmfam=\twelvebm
    \scriptfont\bmfam=\ninebm\scriptscriptfont\bmfam=\sevenbm
  \baselineskip=17.8pt \rm 
  \def\looselineskip{\baselineskip=18.5pt plus 1.8pt}%
  \def\tightlineskip{\baselineskip=16.5pt}%
  \def\verytightlineskip{\baselineskip=15pt}%
  \let\big=\twelvebig\let\Big=\twelveBig\let\bigg=\twelvebigg
  \let\Bigg=\twelveBigg  }
\def\elevenpoint{
  \textfont0=\elevenrm \scriptfont0=\ninerm \scriptscriptfont0=\sixrm
  \def\rm{\fam0 \elevenrm}%
  \textfont1=\eleveni \scriptfont1=\ninei \scriptscriptfont1=\sixi
  \textfont2=\elevensy \scriptfont2=\ninesy \scriptfont2=\sixsy 
  \textfont3=\elevenex \scriptfont3=\elevenex \scriptfont3=\elevenex
  \def\it{\fam\itfam\elevenit}\textfont\itfam=\elevenit
  \def\sl{\fam\slfam\elevensl}\textfont\slfam=\elevensl
  \def\bf{\fam\bffam\elevenbf}\textfont\bffam=\elevenbf
  \def\tt{\fam\ttfam\eleventt}\textfont\ttfam=\eleventt
  \def\sc{\fam\scfam\elevensc}\textfont\scfam=\elevensc
  \def\am{\fam\amfam\elevenam}\textfont\amfam=\elevenam
    \scriptfont\amfam=\nineam\scriptscriptfont\amfam=\sixam
  \def\bm{\fam\bmfam\elevenbm}\textfont\bmfam=\elevenbm
    \scriptfont\bmfam=\ninebm\scriptscriptfont\bmfam=\sixbm
  \baselineskip=15.1pt \rm
  \def\looselineskip{\baselineskip=16pt plus 1.5pt}%
  \def\tightlineskip{\baselineskip=14pt}%
  \def\verytightlineskip{\baselineskip=13pt}%
  \let\big=\elevenbig\let\Big=\elevenBig\let\bigg=\elevenbigg
  \let\Bigg=\elevenBigg  }
\def\tenpoint{
  \textfont0=\tenrm \scriptfont0=\eightrm \scriptscriptfont0=\fiverm
  \def\rm{\fam0 \tenrm}%
  \textfont1=\teni \scriptfont1=\eighti \scriptscriptfont1=\fivei
  \textfont2=\tensy \scriptfont2=\eightsy \scriptfont2=\fivesy 
  \textfont3=\tenex \scriptfont3=\tenex \scriptfont3=\tenex
  \def\it{\fam\itfam\tenit}\textfont\itfam=\tenit
  \def\sl{\fam\slfam\tensl}\textfont\slfam=\tensl
  \def\bf{\fam\bffam\tenbf}\textfont\bffam=\tenbf
  \def\tt{\fam\ttfam\tentt}\textfont\ttfam=\tentt
  \def\sc{\fam\scfam\tensc}\textfont\scfam=\tensc
  \def\am{\fam\amfam\tenam}\textfont\amfam=\tenam
    \scriptfont\amfam=\eightam \scriptscriptfont\amfam=\fiveam
  \def\bm{\fam\bmfam\tenbm}\textfont\bmfam=\tenbm
    \scriptfont\bmfam=\eightbm \scriptscriptfont\bmfam=\fivebm
  \baselineskip=14pt\rm
  \def\looselineskip{\baselineskip=14.8pt plus1.5pt}
  \def\tightlineskip{\baselineskip=12.6pt}%
  \def\verytightlineskip{\baselineskip=13pt}%
  \let\big=\tenbig\let\Big=\tenBig\let\bigg=\tenbigg\let\Bigg=\tenBigg  }
\def\ninepoint{
  \textfont0=\ninerm \scriptfont0=\sevenrm \scriptscriptfont0=\fourrm
  \def\rm{\fam0 \ninerm}%
  \textfont1=\ninei \scriptfont1=\seveni \scriptscriptfont1=\fouri
  \textfont2=\ninesy \scriptfont2=\sevensy \scriptfont2=\foursy 
  \textfont3=\nineex \scriptfont3=\nineex \scriptfont3=\nineex
  \def\it{\fam\itfam\nineit}\textfont\itfam=\nineit
  \def\sl{\fam\slfam\ninesl}\textfont\slfam=\ninesl
  \def\bf{\fam\bffam\ninebf}\textfont\bffam=\ninebf
  \def\tt{\fam\ttfam\ninett}\textfont\ttfam=\ninett
  \def\sc{\fam\scfam\ninesc}\textfont\scfam=\ninesc
  \def\am{\fam\amfam\nineam}\textfont\amfam=\nineam
    \scriptfont\amfam=\nineam\scriptscriptfont\amfam=\fouram
  \def\bm{\fam\bmfam\ninebm}\textfont\bmfam=\ninebm
    \scriptfont\bmfam=\ninebm\scriptscriptfont\bmfam=\fourbm
  \baselineskip=12.6pt\rm
  \def\tightlineskip{\baselineskip=11.5pt}
  \let\big=\ninebig\let\Big=\nineBig\let\bigg=\ninebigg
  \let\Bigg=\nineBigg  }
\def\eightpoint{
  \textfont0=\eightrm \scriptfont0=\fiverm \scriptscriptfont0=\fourrm
  \def\rm{\fam0 \eightrm}%
  \textfont1=\eighti \scriptfont1=\fivei \scriptscriptfont1=\fouri
  \textfont2=\eightsy \scriptfont2=\fivesy \scriptfont2=\foursy 
  \textfont3=\eightex \scriptfont3=\eightex \scriptfont3=\eightex
  \def\it{\fam\itfam\eightit}\textfont\itfam=\eightit
  \def\sl{\fam\slfam\eightsl}\textfont\slfam=\eightsl
  \def\bf{\fam\bffam\eightbf}\textfont\bffam=\eightbf
  \def\tt{\fam\ttfam\eighttt}\textfont\ttfam=\eighttt
  \def\sc{\fam\scfam\eightsc}\textfont\scfam=\eightsc
  \def\am{\fam\amfam\eightam}\textfont\amfam=\eightam
    \scriptfont\amfam=\eightam\scriptscriptfont\amfam=\fouram
  \def\bm{\fam\bmfam\eightbm}\textfont\bmfam=\eightbm
    \scriptfont\bmfam=\eightbm\scriptscriptfont\bmfam=\fourbm
  \baselineskip=11.2pt \rm
  \def\tightlineskip{\baselineskip=10.4pt}
  \let\big=\eightbig\let\Big=\eightBig\let\bigg=\eightbigg
  \let\Bigg=\eightBigg  }

\twelvepoint
\nopagenumbers
\hsize=6in\vsize=8.8in

\parskip=1pt plus 1pt

\newif\ifSpecialhead\Specialheadfalse
\newbox\specialheadbox

\def\specialhead #1\par{\Specialheadtrue\setbox\specialheadbox=\hbox{#1}}
\headline={{\ifSpecialhead\box\specialheadbox\global\Specialheadfalse\else
     \ifnum\pageno<0{\hfill\quad{\twelvebf\folio}}%
     \else\ifnum\pageno<2\hfill
     \else\hfill\twelvepoint\sc\firstmark\quad{\twelvebf\folio}\fi\fi\fi}}

\def\title#1\par{\bigskip{\def\cr{\par\center}\center\fifteenbf #1\par}\medskip}
\def\subtitle#1\par{\centerline{\fifteenrm #1}\medskip}
\def\author#1\par{\medskip{\def\cr{\par\center\twelvesc}\fifteensc\center#1\par}}
\def\center#1\par{\hfil #1\hfil\par}
\def\abstract.#1\par{\message{Abstract.}%
                    \medskip{\narrower\narrower\tenpoint\tightlineskip
                        \noindent{\bf Abstract.}#1\par}\medskip\noindent}
\def\tinyabstract.#1\par{\message{Abstract.}%
                    \medskip{\narrower\narrower\eightpoint\tightlineskip
                        \noindent{\bf Abstract.}#1\par}\medskip\noindent}
\def\bigabstract.#1\par{\message{Abstract.}%
                         \medskip{\narrower\narrower\tightlineskip
                         \noindent{\bf Abstract. }#1\par}\medskip\noindent}
\def\acknowledgement#1\par{\footnote{}{#1}}
\def\sectionskip{\Goodbreak\vskip 25pt plus 15pt minus 5pt}
\def\secnumber{\ifquiet
               \else\ifNoSections
                    \else\sectionsymbol\the\secno\quad\fi\fi}
\def\section#1\par{ \NoSectionsfalse\par\sectionskip\proofdepth=0\claimno=0
 \ifquiet\else\advance\secno by1\fi\toks0={#1}
 \immediate\write16{\ifquiet\else Section \the\secno\space\fi
                    \the\toks0}%
 \mark{\secnumber #1}%
 {\fifteenpoint\bf\noindent\secnumber #1}\nobreak\bigskip\quietoff
 \nobreak\noindent}

\def\QUIET{\QUIETtrue\quiettrue}

\def\quietoff{\ifQUIET\else\quietfalse\fi}
\newif\ifquiet
\newif\ifQUIET
\newif\ifNoSections
\newcount\claimtype
\newcount\secno
\newcount\claimno
\newcount\subclaimno
\newcount\subsubclaimno
\newcount\subsubsubclaimno
\newcount\proofdepth
\def\subclaimnumber{\ifquiet\else\ifcase\subclaimno\or A\or B\or C\or D\or E\or
     F\or G\or H\or I\or J\or K\or L\or M\or N\or O\or P\fi\fi}
\def\subsubclaimnumber{\ifquiet\else\ifcase\subsubclaimno\or i\or ii\or iii\or 
   iv\or v\or vi\or vii\or viii\or ix\or x\or xi\or xii\or xiii\or xiv\fi\fi}
\def\subsubsubclaimnumber{\ifquiet\else\ifcase\subsubsubclaimno\or a\or b\or 
   c\or d\or e\or f\or g\or \or h\or i\or j\or k\or l\or m\or n\or o\fi\fi}
\def\claimtag{\ifquiet\else
  \ifNoSections
    \ifcase\proofdepth\the\claimno%
    \or\the\claimno.\subclaimnumber
    \or\the\claimno.\subclaimnumber.\subsubclaimnumber
    \or\the\claimno.\subclaimnumber.\subsubclaimnumber
                                                .\subsubsubclaimnumber\fi
  \else
    \ifcase\proofdepth\the\secno.\the\claimno
    \or\the\secno.\the\claimno.\subclaimnumber
    \or\the\secno.\the\claimno.\subclaimnumber.\subsubclaimnumber
    \or\the\secno.\the\claimno.\subclaimnumber.\subsubclaimnumber
                                                .\subsubsubclaimnumber\fi\fi\fi}
\secno=0\claimno=0\proofdepth=0\subclaimno=0\subsubclaimno=0\subsubsubclaimno=0
\NoSectionstrue
\newbox\qedbox
\def\claimname{\ifcase\claimtype Theorem\or Lemma\or Claim\or Corollary\or
               Question\or Definition\or Remark\or Conjecture\fi}
\def\preclaimskip{\removelastskip
    \ifcase\claimtype\goodbreak\vskip 8pt plus 4pt minus 2pt
                  \or\goodbreak\vskip 6pt plus 4pt minus 1pt
                  \or\goodbreak\vskip 5pt plus 4pt minus 1pt
                  \or\goodbreak\vskip 8pt plus 4pt minus 2pt
                  \or\vskip 7pt plus 4pt minus 2pt
                  \or\vskip 7pt plus 4pt minus 2pt
                  \or\vskip 7pt plus 4pt minus 2pt
                  \or\goodbreak\vskip 8pt plus 4pt minus 2pt\fi}
\def\postclaimskip{\ifcase\claimtype         \vskip 4pt plus 2pt minus 2pt
                                          \or\vskip 3pt plus 2pt minus 2pt
                                          \or\vskip 2pt plus 2pt minus 1pt
                                          \or\vskip 4pt plus 2pt minus 2pt
                                          \or\vskip 1pt plus 2pt 
                                          \or\vskip 4pt plus 4pt 
                                          \or\vskip 3pt plus 2pt
                                          \or\vskip 4pt plus 2pt minus 2pt\fi}
\def\claimfont{\ifcase\claimtype
                  \sl\or\sl\or\sl\or\sl\or\sl\or\rm\or\rm\or\sl\fi}
\def\advancetag{\ifcase\proofdepth\advance\claimno by1
                               \or\advance\subclaimno by1
                               \or\advance\subsubclaimno by1
                               \or\advance\subsubsubclaimno by1\fi}
\def\sayclaim#1.#2 #3\par{\ifquiet\else\advancetag\fi
    \preclaimskip\setbox1=\hbox{#1}\setbox2=\hbox{#2}%
    \toks0={#1 }
    \immediate\write16{\ifdim\wd1>0pt\the\toks0
                       \else\claimname\space\fi \claimtag.}%
    \vbox{\noindent
    {\bf\ifdim\wd1=0pt \claimname\else #1\fi\ifquiet.\else\ \claimtag{\ifNoSections.\fi}\fi}%
    \enspace{\ifdim\wd2>0pt\sc #2\enspace\fi}%
    {\claimfont #3\par}}\postclaimskip\quietoff}
\def\theorem{\claimtype=0\sayclaim}

\def\corollary{\claimtype=3\sayclaim}

\def\point#1. #2\par{\item{\rm #1.}#2\par}
\def\points#1\cr\par{\medskip\vbox{\let\cr=\point\point#1\par}\par}

\def\prooffont{}
\def\proofsize{}
\def\proofindent{}
\def\proofskip{\badbreak\ifcase\claimtype    \vskip 3pt plus 2pt minus 2pt
                                          \or\vskip 2pt plus 2pt minus 2pt
                                          \or\vskip 1pt plus 2pt minus 1pt
                                          \or\vskip 3pt plus 2pt minus 2pt
                                          \or\vskip 1pt plus 2pt 
                                          \or\vskip 2pt plus 4pt 
                                          \or\vskip 1pt plus 2pt
                                          \or\vskip 3pt plus 2pt minus 2pt\fi}

\def\Goodbreak{\vskip0pt plus.5in\penalty-1000\vskip0pt plus-.5in}
\def\goodbreak{\penalty-500}
\def\badbreak{\penalty500}
\def\Badbreak{\penalty1000}
\def\proof{\message{proof}\removelastskip\Badbreak\proofskip\begingroup
  \advance\proofdepth by1
  \setbox\qedbox=\hbox{\halmos\raise2pt\hbox{\fiverm\claimname}}%
  \prooffont\proofsize\proofindent\noindent{\bf Proof: }}
\def\proofof#1:{\message{proof}\removelastskip\Badbreak\proofskip\begingroup
  \advance\proofdepth by1
  \setbox\qedbox=\hbox{\halmos\raise2pt\hbox{\fiverm#1}}%
  \prooffont\proofsize\proofindent\noindent{\bf Proof of #1: }}
\def\cite[#1]{[{\tenrm{#1}}]\message{[#1]}}
\edef\ref#1{\expandafter\global\expandafter\edef#1{\noexpand\claimtag}}
\newwrite\notes
\openout\notes=\jobname.notes
\long\def\unexpandedwrite#1#2{\def\finwrite{\write#1}%
   {\aftergroup\finwrite\aftergroup{\sanitize#2\endsanity}}}
\def\sanitize{\futurelet\next\sanswitch}
\let\stoken=\space
\def\sanswitch{\ifx\next\endsanity
  \else\ifcat\noexpand\next\stoken\aftergroup\space\let\next=\eat
   \else\ifcat\noexpand\next\bgroup\aftergroup{\let\next=\eat
    \else\ifcat\noexpand\next\egroup\aftergroup}\let\next=\eat
     \else\let\next=\copytoken\fi\fi\fi\fi \next}
\def\eat{\afterassignment\sanitize \let\next= }
\long\def\copytoken#1{\ifcat\noexpand#1\relax\aftergroup\noexpand
  \else\ifcat\noexpand#1\noexpand~\aftergroup\noexpand\fi\fi
  \aftergroup#1\sanitize}
\def\endsanity\endsanity{}

\def\note#1#2{\hbox to2in{\strut#1\quad\dotfill\quad#2}}
\def\boxit#1{\setbox4=\hbox{\kern1pt#1\kern1pt}
  \hbox{\vrule\vbox{\hrule\kern1pt\box4\kern1pt\hrule}\vrule}}
\def\halmos{\hbox{\am\char'3}} 
\def\qed#1\par{\message{.                                }\setbox1=\hbox{#1}%
  \ifdim\wd1>0pt\setbox\qedbox=\hbox{\halmos\raise2pt\hbox{\fiverm #1}}\fi
  \kern5pt\lower 2pt\hbox{\box\qedbox}\proofskip\goodbreak\endgroup}

\def\sectionsymbol{\S}
\def\k{\kappa}
\def\g{\gamma}
\def\a{\alpha}

\def\l{\lambda}

\def\I1{\mathop{\hbox{\sc i}_1}}
\def\w{\omega}
\def\P{{\mathchoice{\hbox{\bm P}}{\hbox{\bm P}}
         {\hbox{\tenbm P}}{\hbox{\sevenbm P}}}}
\def\Q{{\mathchoice{\hbox{\bm Q}}{\hbox{\bm Q}}
         {\hbox{\tenbm Q}}{\hbox{\sevenbm Q}}}}

\def\elesub{\prec}

\def\unifto{\buildrel\lower 7pt\hbox{$\to$}\over\to}

\def\iso{\cong}

\def\cp{\mathop{\rm cp}\nolimits}

\def\ZFC{\hbox{\sc zfc}}
\def\GCH{\hbox{\sc gch}}

\def\plus{^{\scriptscriptstyle +}}

\def\in{\mathrel{\mathchoice{\raise 
1pt\hbox{$\scriptstyle\cal\char'62$}}
         {\raise 1pt\hbox{$\scriptstyle\cal\char'62$}}
         {\raise .5pt\hbox{$\scriptscriptstyle\cal\char'62$}}
         {\hbox{$\scriptscriptstyle\cal\char'62$}}}\penalty700{}}
\def\ni{\mathrel{\mathchoice{\raise 1pt\hbox{$\scriptstyle\cal\char'63$}}
                   {\raise 1pt\hbox{$\scriptstyle\cal\char'63$}}
                   {\raise .5pt\hbox{$\scriptscriptstyle\cal\char'63$}}
                   {\hbox{$\scriptscriptstyle\cal\char'63$}}}\penalty700}
\def\of{\mathrel{\mathchoice{\raise 1pt\hbox{$\scriptstyle\subseteq$}}
                   {\raise 1pt\hbox{$\scriptstyle\subseteq$}}
                   {\raise .5pt\hbox{$\scriptscriptstyle\subseteq$}}
                   {\hbox{$\scriptscriptstyle\subseteq$}}}}
\def\fo{\mathrel{\mathchoice{\raise 1pt\hbox{$\scriptstyle\supseteq$}}
                   {\raise 1pt\hbox{$\scriptstyle\supseteq$}}
                   {\raise .5pt\hbox{$\scriptscriptstyle\supseteq$}}
                   {\hbox{$\scriptscriptstyle\supseteq$}}}}
\def\notin{\mathrel{\mathchoice
  {\raise 1pt\hbox{\rlap{$\scriptstyle\;|$}$\scriptstyle\cal\char'62$}}
  {\raise 1pt\hbox{\rlap{$\scriptstyle\kern2pt 
          |$}$\scriptstyle\cal\char'62$}}
  {\raise .5pt\hbox{\rlap{$\scriptscriptstyle\, |$}$\scriptscriptstyle
      \cal\char'62$}}
  {\hbox{\rlap{$\scriptscriptstyle\, |$}$\scriptscriptstyle
     \cal\char'62$}}}%
  \penalty700}

\def\and{\mathrel{\kern1pt\&\kern1pt}}
\def\iff{\mathrel{\leftrightarrow}}

\def\union{\cup}

\def\[#1]{\left[\vphantom{\bigm|}#1\right]}
\def\<#1>{\langle\,#1\,\rangle}

\def\sat{\models}

\def\image{\mathbin{\hbox{\tt\char'42}}}
\def\restrict{\mathbin{\mathchoice{\hbox{\am\char'26}}{\hbox{\am\char'26}}{\hbox{\eightam\char'26}}{\hbox{\sixam\char'26}}}}

\def\st{\mid}
\def\seq<#1>{{\def\st{\mid\penalty650}\left<\,#1\,\right>}}

\def\set#1{\{\,#1\,\}}

\def\th{{\hbox{\fiverm th}}}

\def\lttheta{{\raise 1pt\hbox{$\scriptstyle<$}\theta}}

\def\I1{\mathop{\hbox{\sc i}_1}}

\def\lteg{{{\scriptstyle\leq}\g}}

\def\ltek{{\scriptstyle\leq}\k}

\def\j{\mathop{\hbox{\bf j}}}
\def\WA{\hbox{\sc wa}}
\def\VHOD{\hbox{\sc v=hod}}
\def\I{\hbox{\sc I}}
\def\ltew{{\scriptstyle\leq}\w}
\def\LST{\{{\in\!}\}}
\def\LWA{\{{\in\!},\j\}}
\QUIET
\center [submitted to the Archive for Mathematical Logic]

\title The Wholeness Axioms and V=HOD

\author Joel David Hamkins\cr
            Kobe University and\cr
            The City University of New York\cr
	{\tentt http://www.math.csi.cuny.edu/$\sim$hamkins}\cr

\acknowledgement My research has been supported in part by a grant from the PSC-CUNY Research Foundation and a fellowship from the Japan Society for the Promotion of Science. And I would like to thank my gracious hosts at Kobe University for their generous hospitality.

\abstract. If the Wholeness Axiom $WA_0$ is itself consistent, then it is consistent with \VHOD. A consequence of the proof is that the various Wholeness Axioms are not all equivalent. Additionally, the theory $\ZFC+\WA_0$ is finitely axiomatizable. 

The Wholeness Axioms, proposed by Paul Corazza, occupy a high place in the upper stratosphere of the large cardinal hierarchy. They are intended as weakenings of the famous inconsistent assertion that there is a nontrivial elementary embedding from the universe to itself, weakenings which, one hopes, are substantial enough to avoid inconsistency but slight enough for them to remain very strong. 
The Wholeness Axioms are formalized in the language $\LWA$, augmenting the usual language of set theory $\LST$ with an additional unary function symbol $\j$ to represent the embedding. The base theory \ZFC\ is expressed only in the smaller language $\LST$. Corazza's original proposal, which I will denote by $\WA_0$, asserts that $\j$ is a nontrivial amenable elementary embedding from the universe to itself. Elementarity is expressed by the scheme $\varphi(x)\iff\varphi(\j(x))$, where $\varphi$ runs through the formulas of the usual language of set theory; nontriviality is expressed by the sentence $\exists x \j(x)\not=x$; and amenability is simply the assertion that $\j\restrict A$ is a set for every set $A$. One can easily see that amenability in this case is equivalent to the assertion that the Separation Axiom holds for $\Sigma_0$ formulae in the language $\LWA$. Using this idea and hankering for a stronger assumption, Corazza finally settled on the version of the Wholeness Axiom that I will here denote by $\WA_\infty$, which asserts in addition that the full Separation Axiom holds in the language $\LWA$. 

As I hope my notation suggests, these two axioms are the important endpoints of a natural hierarchy of axioms $\WA_0$, $\WA_1$, $\WA_2,\ldots,\WA_\infty$, which I will refer to collectively as the Wholeness Axioms. Specifically, the Wholeness Axiom $\WA_n$, where $n$ is amongst $0,1,\ldots,\infty$, consists of the following:

\points 1. (Elementarity) All instances of $\varphi(x)\iff\varphi(\j(x))$ for $\varphi$ in the language $\LST$.\cr
           2. (Separation) All instances of the Separation Axiom for $\Sigma_n$ formulae in the full language $\LWA$.\cr
           3. (Nontriviality) The axiom $\exists x\,\j(x)\not=x$.\cr

Those familiar with the Kunen \cite[Kun71] inconsistency may fear to tread so near calamity! But the attentive reader will notice that what is missing from the Wholeness Axiom schemes, and what figures prominantly in Kunen's proof, are the instances of the Replacement Axiom in the full language with $\j$. And this difference is important. Kunen uses the Replacement Axiom in the full language in order to know that the critical sequence $\set{\k_n\st n\in\w}$, defined as usual by $\k_0=\k=\cp(j)$ and $\k_{n+1}=j(\k_n)$, is a set. But in the language without the symbol $\j$, there is no way to define the critical sequence, and so under the Wholeness Axiom Kunen's argument cannot be carried out. 

There is occasionally the suggestion in the literature that we should take the Kunen \cite[Kun71] inconsistency only to rule out definable embeddings, that is, embeddings $j:V\to V$ which as a class are definable from parameters by a formula, and that we should understand Kunen's theorem as a scheme, one theorem for each possible definition of $j$. But to do so weakens the force of Kunen's result. Quite aside from Kunen's argument, it is relatively easy to rule out definable embeddings from the universe to itself (and \cite[Suz$\infty$] makes this point plain). The argument is simply this: given a formula $\varphi(x,y,\vec z)$ purportedly defining a nontrivial embedding $j:x\mapsto y$ from $V$ to $V$, possibly with parameters $\vec z$, let $\k$ be least such that for some parameters $\vec a$ the formula $\varphi(x,y,\vec a)$ defines an elementary embedding from $V$ to $V$ with critical point $\k$ (this is seen to be first-order expressible using a $\Sigma_1$ truth predicate and the fact that any cofinal $\Sigma_1$-elementary embedding is fully elementary). Applying now the embedding $j$, we see that $j(\k)$ is also the least possible critical point of an elementary embedding of $V$ defined by $\varphi$, and this is a contradiction since $\k<j(\k)$. The point of this simple argument is that one doesn't need Kunen's argument to rule out definable embeddings; indeed, the argument shows that one does not even need the Axiom of Choice to do so. Rather, we should view the Kunen inconsistency as taking place in the language with the extra unary function symbol $\j$, and see that what he has really proved is that the axioms of \ZFC\ in the language $\set{{\in},\j}$ are inconsistent with the theory asserting that $\j$ represents a nontrivial elementary embedding from $V$ to $V$, as in 1 and 3 above. Taken this way, Kunen's result does not assume that the embedding is definable and leads naturally to the Wholeness Axioms as straightforward weakenings of an inconsistent theory.

Since as I have said $\WA_0$ already implies the amenability of the embedding, every model of one of the Wholeness Axioms has the form $\<V,{\in},j>$, where $\<V,{\in}>$ is a model of $\ZFC$ and $j:V\to V$ is a nontrivial amenable elementary embedding. Kunen's inconsistency proof is easily modified to show that under the Wholeness Axioms, the critical sequence $\set{\k_n\st n\in\w}$ is unbounded in the ordinals of $V$. And as Corazza has noted, the embedding $j$ of course cannot be definable in $V$. Since $j\restrict V_\k$ is the identity function, it follows that $V_\k\elesub V_{\k_1}$, and by iteratively applying $j$ to this fact one easily concludes that $$V_{\k_0}\elesub V_{\k_1}\elesub V_{\k_2}\elesub\cdots\elesub V.$$
So under the Wholeness Axioms the universe $V$ is, externally, the union of this countable elementary chain of end-extensions. 

For an upper bound on the consistency strength of the Wholeness Axioms, one can easily show that if $\I_3$ holds, that is, if there is a nontrivial elementary embedding $j:V_\l\to V_\l$, then $\<V_\l,{\in},j>$ is a model of the Wholeness Axiom $\WA_\infty$ (see \cite[Cor$\infty$a]). Conversely, for a lower bound it is easy to see that if $\k$ is the critical point of the embedding $j$ arising in the Wholeness Axiom $\WA_0$, then $\k$ is supercompact, extendible, etc. And using the same kind of argument, Corazza \cite[Cor$\infty$a] has shown that such a $\k$ is super-$n$-huge for every natural number $n$.

Let me now turn to the main focus of this paper, namely, the question of whether the Wholeness Axioms are consistent with \VHOD. Corazza \cite[Cor$\infty$b] has proved:

\theorem.(Corazza) If there is an $\I_1$ embedding $j:V_{\l+1}\to V_{\l+1}$, then the Wholeness Axiom $\WA_\infty$ is consistent with \VHOD. Hence also the Wholeness Axiom $\WA_0$ is consistent with \VHOD. 

Since the existence of an $\I_1$ embedding is strictly stronger than $\WA_\infty$ in consistency strength, one naturally hopes to improve the hypothesis of the theorem to a pure relative consistency result, and Corazza conjectures as much in \cite[Cor$\infty$b]. That is, for the optimal hypothesis, one would want to assume only the consistency of the Wholeness Axiom itself and make the same conclusion. And in the case of the Wholeness Axiom $\WA_0$, that is what I achieve in this paper:

\theorem Main Theorem. If the Wholeness Axiom $\WA_0$ is itself consistent, then it is consistent with \VHOD. 

\proof Of course, I am speaking here about consistency over the base theory \ZFC. For simplicity, let me first obtain a model of the Wholeness Axiom $\WA_0$ in which the \GCH\ holds. After this, I will show fully how to obtain \VHOD. So suppose that $V$ is a model of the Wholeness Axiom $\WA_0$, with embedding $j:V\to V$. Let me remark that this argument will not take place inside $V$; rather, I will freely make use of second-order manipulations of $V$ and $j$ in order to construct the desired model. As usual, let $\<\k_n\st n\in\w>$ be the critical sequence of $j$, so that $\k=\k_0=\cp(j)$ and $\k_{n+1}=j(\k_n)$. I mentioned earlier that 
$$V_{\k_0}\elesub V_{\k_1}\elesub V_{\k_2}\elesub\cdots\elesub V,$$
and moreover $V$ is the union of this elementary chain. Let $\P_\k$ be the usual reverse Easton $\k$-iteration in $V$ which forces the \GCH\ up to $\k$. Thus, at any stage $\g$ which is a cardinal in $V[G_\g]$, the forcing $\Q_\g$ adds a Cohen subset to $\g\plus$, and we take direct limits at inaccessible stages and otherwise inverse limits. The stage $\g$ forcing, of course, ensures that $2^\g=\g\plus$, and for any $\g<\k$ we may factor the entire iteration as $\P_\k\iso\P_\g*\P_{\g,\k}$, where $\P_{\g,\k}$ is $\lteg$-directed closed in $V^{\P_\g}$. So $\P_\k$ will naturally force the \GCH\ up to $\k$. Rather than taking a formal direct limit at stage $\k$ using equivalence classes of threads, etc., let me more simply regard $\P_\k$ as a subset of $V_\k$ by identifying each thread with it's initial point. By this device, names for elements of rank less than $\k$ themselves have rank less than $\k$. Henceforth denoting $\k$ by $\k_0$, suppose now that $G_{\k_0}\of\P_{\k_0}$ is $V$-generic, and consider the forcing $j(\P_{\k_0})=\P_{\k_1}$. We may factor this longer iteration as $\P_{\k_1}\iso\P_{\k_0}*\P_{\k_0,\k_1}$, where $\P_{\k_0,\k_1}$ is $\ltek_0$-directed closed in $V[G_{\k_0}]$. Forcing to add a $V[G_{\k_0}]$-generic $G_{\k_0,\k_1}\of\P_{\k_0,\k_1}$ provides a $V$-generic filter $G_{\k_1}=G_{\k_0}*G_{\k_0,\k_1}$. And since $j\image G_{\k_0}=G_{\k_0}\of G_{\k_1}$ we may lift the embedding to $j:V[G_{\k_0}]\to V[G_{\k_1}]$ with $j(G_{\k_0})=G_{\k_1}$. In particular, observe that $V_{\k_0}[G_{\k_0}]\elesub V_{\k_1}[G_{\k_1}]$. Now consider the forcing $j(\P_{\k_1})\iso\P_{\k_0}*\P_{\k_0,\k_2}\iso\P_{\k_1}*\P_{\k_1,\k_2}$. We already have the generic filter $G_{\k_1}$ for the first $\k_1$ many stages of forcing. Furthermore, the remaining forcing $\P_{\k_1,\k_2}$ is $\ltek_1$-directed closed in $V[G_{\k_1}]$. Furthermore, $j\image G_{\k_0,\k_1}$ is a directed subset of this poset in $V[G_{\k_1}]$ of size $\k_1$, so by the directed closure of $\P_{\k_1,\k_2}$ in $V[G_{\k_1}]$ there is a (master) condition $q$ below every element of $j\image G_{\k_0,\k_1}$. Force below this condition to add a generic $G_{\k_1,\k_2}$. This provides a generic $G_{\k_2}\iso G_{\k_1}*G_{\k_1,\k_2}$, and because of the master condition we know $j\image G_{\k_1}\of G_{\k_2}$. Continuing in this manner, we may obtain generic objects $G_{\k_n}$ for every natural number $n$ with the properties that $G_{\k_n}\of P_{\k_n}$ is $V$-generic and $j\image G_{\k_n}\of G_{\k_{n+1}}$.

This means that we may always lift the embedding to $j:V[G_{\k_n}]\to V[G_{\k_{n+1}}]$, with $j(G_{\k_n})=G_{\k_{n+1}}$. (While this embedding is defined in $V[G_{\k_{n+1}}]$, the domain is only $V[G_{\k_n}]$.) Since as I have mentioned earlier $V_{\k_0}[G_{\k_0}]\elesub V_{\k_1}[G_{\k_1}]$, we may now iteratively apply $j$ to conclude that $V_{\k_n}[G_{\k_n}]\elesub V_{\k_{n+1}}[G_{\k_{n+1}}]$. In summary, we have constructed an elementary chain of end extensions: 
$$V_{\k_0}[G_{\k_0}]\elesub V_{\k_1}[G_{\k_1}]\elesub V_{\k_2}[G_{\k_2}]\elesub\cdots$$
Now I will make the key step of the proof. Let 
$$\tilde V=\union_{n\in\w}V_{\k_n}[G_{\k_n}].$$ Since $\tilde V$ is the union of this elementary chain, the theory of $\tilde V$ is the same as the theory of each $V_{\k_n}[G_{\k_n}]$. In particular, $\tilde V$ is a model of \ZFC\ in which the \GCH\ holds, since by construction $G_\k$ ensures the \GCH\ up to $\k$ in $V_\k[G_\k]$. Also, we have already extended $j$ to each $j:V[G_{\k_n}]\to V[G_{\k_{n+1}}]$, so we have a map defined $j:\tilde V\to \tilde V$. Let me now argue that this map is elementary. If $\tilde V\sat\varphi(x)$, then by the elementary chain $V_{\k_n}[G_{\k_n}]\sat\varphi(x)$ for sufficiently large $\k_n$; thus, by applying $j$ we conclude that $V_{\k_{n+1}}[G_{\k_{n+1}}]\sat\varphi(j(x))$ and so by the elementary chain again we see that $\tilde V\sat\varphi(j(x))$. Further, since $j\restrict V[G_{\k_n}]$ was defined in $V[G_{\k_{n+1}}]$, it follows that $j\restrict V_{\k_n}[G_{\k_n}]\in\tilde V$ and so the embedding is amenable. Thus, $\<\tilde V,\in,j>$ is a model of $\WA_0$ and the \GCH. 

Now I will explain fully how to obtain \VHOD. Suppose that $V$ is a model of the Wholeness Axiom $\WA_0$ in which the \GCH\ also holds, and that $j:V\to V$ is the witnessing embedding, with critical sequence as above. The idea is to force \VHOD\ by coding every set into the continuum function. That is, we will make the \GCH\ hold or fail at each $\aleph_\a$ so as to code one bit of information. Let $\P_\k$ be the reverse Easton $\k$-iteration which uses this technique to code every set of rank less than $\k$ into the continuum function below $\k$ (the same forcing is used with an $\I_1$ embedding in \cite[Cor$\infty$b]). Thus, at an inaccessible stage of forcing $\g<\k$, we select, using a fixed well-ordering of the names in $V_\k$, a set $A\of\g$ coding a relation $E$ on $\g$ such that $\<V_\g[G_\g],{\in}>\iso\<\g,E>$. Next, we code $A$ and hence $E$ and hence $V_\g[G_\g]$ into the continuum function above $\g$ by forcing the continuum hypothesis to hold or fail at $\aleph_{\g+\a+1}$ depending on whether $\a\in A$ or not. For the next nontrivial stage of forcing, we wait until the next inaccessible cardinal, which is beyond all this coding, so that the various stages of coding do not interfere with each other. Since unboundedly often we make sets definable from the continuum function, this iteration will force \VHOD\ in $V_\k^{\P_\k}$. Also, the iteration has the same closure properties as the previous iteration by which we obtained the \GCH. By the identical lifting arguments as I just gave in that case, therefore, we may again construct the model $\tilde V=\union_{n\in\w}V_{\k_n}[G_{\k_n}]$ as the union of an elementary chain of end-extensions. The model $\tilde V$ again has the same theory as $V_\k[G_\k]$. In particular, it is a model of \ZFC\ in which \VHOD\ holds. And again the embedding lifts to $j:\tilde V\to\tilde V$ in such a way so as to witness the Wholeness Axiom $\WA_0$. So the proof is complete.\qed

The technique shows that if the Wholeness Axiom $\WA_0$ holds, where the embedding has critical point $\k$, and $\P_\k$ is a reverse Easton $\k$-iteration with $\lteg$-directed closed forcing at every stage $\g$, then any fact true in $V_\k^{\P_\k}$ is consistent with $\WA_0$. To see this, one simply constructs the generics $G_{\k_n}\of\P_{\k_n}$ as above and lifts the Wholeness embedding to $\tilde V=\union_nV_{\k_n}[G_{\k_n}]$. The lifted embedding witnesses $\WA_0$ in $\tilde V$ and the theory of $\tilde V$ is the same as the theory of $V_\k[G_\k]$ because of the elementary chain. 

The question remains open whether the relative consistency results hold at each level, that is, whether the consistency of $\WA_n$ implies the consistency of $\WA_n$ with $\VHOD$ for $1\leq n\leq\infty$. Though my proof above does not appear to settle this question, I can extract from this failure the consequence that the various Wholeness Axioms are not all equivalent:

\corollary. If consistent, the Wholeness axiom $\WA_0$ is not logically equivalent to any other $\WA_n$.

\proof What I will show is that we can arrange that the model $\tilde V$ in the proof of the Main Theorem satisfies $\WA_0$ but not $\WA_1$. Since the other $\WA_n$ all imply $\WA_1$, this will establish the corollary. 

The first step is to fix a real $u\of\w$ which is not in $V$. For example, perhaps $u$ is generic over $V$. More deviously, perhaps $u$ identifies indiscernibles for $V$ or, in a truly naughty case, perhaps $u$ collapses every cardinal of $V$ to $\w$. The point is that while we will not add the naughty real $u$ as an element to $\tilde V$, we will ensure that as a subset of $\w$ it satisfies a $\Sigma_1$ definition there in the language $\LWA$, thereby violating $\Sigma_1$-Separation. 

In the proof of the Main Theorem we must periodically force over $V$ to construct the generics $G_{\k_n,\k_{n+1}}$ below some master condition $q$. This is a kind of free step in the construction, and we will take advantage of this freedom to code one additional bit of the real $u$. Specifically, when we do this forcing there is some first nontrivial stage at which we add a Cohen set (recall that the iteration proceeds by periodically adding a lot of Cohen sets to cardinals in such a way so as to code information into the continuum function). When choosing the generic $G_{\k_n,\k_{n+1}}$, let the first digit of the first such set we add be determined by the $n^\th$ digit of the fixed naughty real $u$. Thus, at the first stage after $\sup j\image\k_n$ at which we add a Cohen set, the first bit of the generic set is the same as $u_n$, the $n^\th$ bit of $u$. Since except for this modification the construction is the same as before, the resulting model $\tilde V$ will be a model of $\WA_0$, just as before.

It remains to show that $\tilde V$ is not a model of $\WA_1$. The first step is to realize that each of the forcing notions $\P_{\k_n}$ is $\ltew$-closed, and so the reals of $\tilde V$ are the same as the reals of $V$. In particular, $u$ is not in $\tilde V$. The second step is to observe that the digits of $u$ may be simply read off from the generics $G_{\k_n}$. Namely, $u_n$ is precisely the first bit of the first Cohen set added after stage $j\image\k_n$ in $G_{\k_{n+1}}$. So $u_n$ is $\Delta_0$-definable in $\tilde V$ from $G_{\k_{n+1}}$ and $\k_n$ (using the function symbol $\j$ to express $j\image\k_n$). Further, each $G_{\k_n}$ is $\Sigma_1$-definable in $\tilde V$ from $G_{\k_0}$, since $G_{\k_n}$ is the last element of the sequence $\seq<X_i\st i\leq n>$ for which $X_0=G_{\k_0}$ and $X_{i+1}=j(X_i)$. Thus, $u$ is a $\Sigma_1$-definable subset (or more properly sub-class) of $\w$ in $\tilde V$, with the parameter $G_{\k_0}$. Since $u$ is not in $\tilde V$, we conclude that $\WA_1$ fails there.\qed

Perhaps the construction of the Main Theorem could be modified in a different way, avoiding all such naughty reals $u$, so as to ensure $\WA_1$ or more generally $\WA_n$ in $\tilde V$, given $\WA_n$ in $V$. This would answer the question I asked earlier about whether the relative consistency result for $\VHOD$ could be proved at the other levels of the Wholeness Axiom hierarchy. 

Let me conclude this paper with the observation that the Wholeness Axiom scheme, at least in the case of $\WA_0$, is actually expressible in a single sentence. 

\theorem. The Wholeness Axiom $\WA_0$ is finitely axiomatizable. Indeed, $\ZFC+\WA_0$ is finitely axiomatizable.

\proof I will use the easily established fact that an elementary embedding is the same thing as a $\Sigma_1$ cofinal embedding, that is, that any map $j:V\to V$ that preserves $\Sigma_1$ truth and has the property that for every $a$ there is a $b$ such that $a\of j(b)$ is fully elementary. This can be proved by a simple induction on formulas. Using now a $\Sigma_1$ truth predicate, one expresses the full elementarity of $j$ by the single assertion: ``$\j$ is a $\Sigma_1$ cofinal map from the universe to itself''. After this, given the full elementarity of $j$, the elementary chain construction shows that $V_\k$ has the same theory as $V$, and so asserting \ZFC\ in $V$ is equivalent to asserting \ZFC\ in $V_\k$, a set, and this is expressible by a single formula. Nontriviality is expressed by the formula $\exists x\, \j(x)\not=x$. Finally, the amenability of $j$ is simply the assertion that for every set $a$ there is a set $b$ which is equal to $j\restrict a$. Thus, $\WA_0$ is finitely axiomatizable.\qed

The question remains open whether the other Wholeness Axioms are finitely axiomatizable.

\section Bibliography

\nopagenumbers
\parindent=0pt
\newbox\Article
\newbox\Journal
\newbox\Author
\newbox\Vol
\newbox\No
\newbox\Year
\newbox\Page
\newbox\Book
\newbox\Publisher
\newbox\Pubaddr
\newbox\Key
\newbox\Editor
\newbox\Comment
\newbox\Note
\def\entry#1#2\par{\item{#1\quad}\hskip-1.1em#2\par}
\def\article#1{\setbox\Article=\hbox{\sl #1, }}
\def\journal#1{\setbox\Journal=\hbox{\rm #1 }}
\def\author#1{\setbox\Author=\hbox{\sc #1, }}
\def\vol#1{\setbox\Vol=\hbox{\bf #1 }}
\def\no#1{\setbox\No=\hbox{no. #1 }}
\def\year#1{\setbox\Year=\hbox{\rm({\oldstyle #1}) }}
\def\page#1{\setbox\Page=\hbox{\rm p. #1 }}
\def\book#1{\setbox\Book=\hbox{\it #1, }}
\def\publisher#1{\setbox\Publisher=\hbox{\rm #1, }}
\def\pubaddr#1{\setbox\Pubaddr=\hbox{\rm #1, }}
\def\key#1{\setbox\Key=\hbox{#1}}
\def\editor#1{\setbox\Editor=\hbox{\rm(#1, Ed.) }}
\def\comment#1{\setbox\Comment=\hbox{\rm #1}}
\def\note#1{\setbox\Note=\hbox{\rm #1 }}
\def\ref#1\par{\smallskip{#1
  \entry{\ifhbox\Key\unhbox\Key\else[\ ]\fi}%
  \unhbox\Author\unhbox\Note
  \ifhbox\Book \unhbox\Book\unhbox\Publisher\unhbox\Pubaddr
               \unhbox\Editor\unhbox\Page\unhbox\Year\unhbox\Comment
  \else \unhbox\Article\unhbox\Journal\unhbox\Vol\unhbox\No\unhbox\Editor
        \unhbox\Page\unhbox\Year\unhbox\Comment\fi\par}}

\tenpoint\tightlineskip

\ref
\author{Paul Corazza}
\article{The Wholeness Axiom and Laver sequences}
\journal{to appear in the Annals of Pure and Applied Logic}
\key{[Cor$\infty$a]}

\ref
\author{Paul Corazza}
\article{The Wholeness Axiom and \VHOD}
\journal{to appear in the Archive for Mathematical Logic}
\key{[Cor$\infty$b]}

\ref
\author{Kenneth Kunen}
\article{Elementary embeddings and infinite combinatorics}
\journal{Journal of Symbolic Logic}
\year{1971}
\vol{36}
\page{407-413}
\key{[Kun71]}

\ref
\author{Akira Suzuki}
\article{No elementary embedding from $V$ to $V$ is definable}
\journal{to appear in the Journal of Symbolic Logic}
\key{[Suz$\infty$]}

\bye